\documentclass[a4paper, 10pt]{article}
\author{Femke Douma}
\title{Spherical Averages on Regular and Semiregular Graphs}
\date{}
\usepackage[dvips]{graphicx}
\usepackage{amsmath}
\usepackage{amssymb}
\usepackage[all]{xy}
\newtheorem{thm}{Theorem}

\newtheorem{lem}{Lemma}
\newtheorem{cor}[lem]{Corollary}

\newtheorem{rmk}[lem]{Remark}

\newenvironment{pf}{\noindent \textsc{Proof} \rm}

\begin{document}

\maketitle

\begin{abstract}
In 1966, P.~G\"unther proved the following result: Given a continuous function $f$ on a compact surface $M$ of constant curvature $-1$ and its periodic lift $\tilde{f}$ to the universal covering, the hyperbolic plane, then the averages of the lift $\tilde{f}$ over increasing spheres converge to the average of the function $f$ over the surface $M$. In this article, we prove similar results for functions on the vertices and edges of regular and semiregular graphs, with special emphasis on the convergence rate. We also consider averages over more general sets like arcs, tubes and horocycles.
\end{abstract}

\section{Introduction and Results}\label{intro}

Let $G$ be a graph with edge set $E$ and vertex set $V$. We require that the graph is finite and connected. It may sometimes have loops and/or multiple edges; if these are not allowed the graph is called simple. Let $d(v)$ be the vertex degree of $v\in V$, where we note that a loop at vertex $v$ contributes $2$ to its degree. We shall only need the edge degree $d'(e)$ for simple graphs $G$, where it is defined as the number of edges meeting $e$ in either of its endpoints. Note that this is equivalent to the vertex degree of $e$ in the line graph $L(G)$ of $G$. Denote by $\delta(v,w)$ the combinatorial distance between vertices $v,w\in V$.

The universal cover $\widetilde{G}$ of any graph $G$ is a tree with vertex set $\widetilde{V}$ and edge set $\widetilde{E}$. It can be constructed as follows: choose a root vertex in $G$, let the vertex set of $\widetilde{G}$ be the set of non-backtracking paths in $G$ starting at the root vertex, and define two such paths to be neighbours if they differ by exactly one edge at the end (see \cite{ow}, p833). The projection map $\pi:\widetilde{G}\rightarrow G$ then associates a vertex $\tilde{v}\in\widetilde{V}$ with the terminal vertex in $V$ of the path it represents, and an edge in $\widetilde{E}$ from $\tilde{v}$ to $\tilde{w}$ is mapped to the edge in $E$ by which the paths $\tilde{v}$ and $\tilde{w}$ differ. If we are now given a real function $f$ on the vertices $V$ (or edges $E$) of $G$, we can lift it uniquely to a function $\tilde{f}$ on the vertices (edges) of the universal cover via $\tilde{f}=f\circ\pi$. 

\begin{figure}[ht]
\centering
\includegraphics[height=7cm]{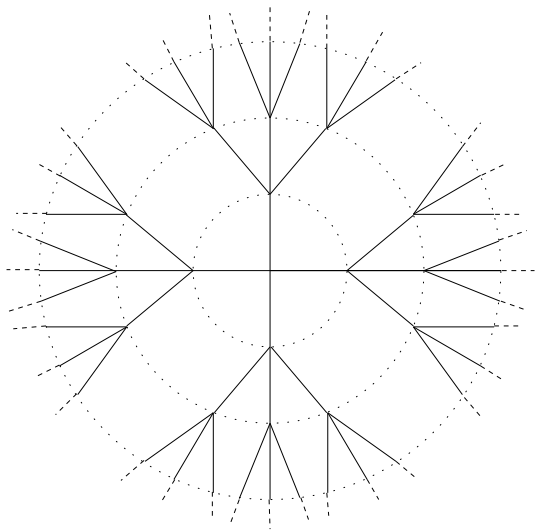}
\caption{Vertex spheres on a regular tree of degree 4}
\end{figure}

We define a vertex sphere  $S_r(v_0) = \{ v \in \widetilde{V} : \delta(v,v_0) = r \}$ on $\widetilde{G}$. The edge sphere is defined analogously as $$S'_r(v_0) = \{ e = \{x,y\} \in \widetilde{E} : \min \{\delta(x,v_0), \delta(y,v_0) \} = r \}$$
The \emph{spherical mean} of the function $f$ on the vertices or the edges is now defined as
$$M_{r,v_0}(f)=\frac{1}{|S_r(v_0)|} \sum_{x\in S_r(v_0)} \tilde{f}(x)$$
where we replace $S_r$ by $S_r'$ in the edge case.

The purpose of this paper is to study the asymptotic behaviour of these averages as $r\rightarrow\infty$ for particular types of graphs, before looking at some other averages. We use an analogy between hyperbolic surfaces and regular graphs which has been studied by various authors (see e.g.~the preface of \cite{figa}). The problem studied by G\"unther \cite{gunther} in the continuous case translates to the case of the regular graph, where we obtain the following result:

\begin{thm}\label{reg v}
Let $G$ be a nonbipartite regular connected graph of degree $d(v)=q+1\geq3$ and $f:V\rightarrow\mathbb{R}$ a function on its vertices. Then we have for any basepoint $v_0\in \widetilde{V}$
$$\Big|M_{r,v_0}(f) - \frac{1}{|V|}\sum_{v\in V} f(v)\Big| \leq  C_G ||f||_2 \beta^r $$
Here $C_G$ is a constant depending on $G$ but independent of $v_0$, and $\beta\in[q^{-1/2},1)$ with lowest value $\beta=q^{-1/2}$ for Ramanujan graphs.
\end{thm}
Obviously, this implies that
$$\lim_{r\rightarrow\infty} M_{r,v_0}(f)= \frac{1}{|V|}\sum_{v\in V} f(v)$$
for any $v_0\in\widetilde{V}$. The norm $||f||_2$ comes from the inner product $\langle f,g \rangle =\sum_{v\in V}f(v)g(v)$. We exclude bipartite graphs in this theorem because spheres of even and odd radii can have different limiting behaviour in this case - see Section \ref{other} for details. We shall see in the proof that $\beta$ depends on the spectral gap - recall that Ramanujan graphs are graphs with a large spectral gap (see e.g.~\cite{valette} or \cite{lubotzky}).

\begin{rmk}\label{arc}
 The proof of Theorem \ref{reg v} can be modified for averages over increasing circular arcs on $\widetilde{G}$. A straightforward argument then allows us to use the result for arcs to prove convergence on increasing subsets of \emph{horocycles} as defined by Cartier in \cite{cartier}. See Section \ref{other} for more details.
\end{rmk}

The next result is concerned with functions defined on the edges of a regular graph. 
\begin{thm}\label{reg e}
Let $G$ be a regular connected simple graph with $d'(e)=2q\geq4$ and let $f:E\rightarrow\mathbb{R}$ be a function on its edges. Then we have for any basepoint $v_0\in \widetilde{V}$
$$\Big|M_{r,v_0}(f) - \frac{1}{|E|}\sum_{e\in E} f(e)\Big| \leq  C_G ||f||_2 \beta^r $$
Here $C_G$ is a constant depending on $G$ but independent of $v_0$, and $\beta\in[q^{-1/2},1)$.
\end{thm}

Note that it is not possible to use Theorem \ref{reg v} to prove Theorem \ref{reg e} (and \ref{semihom} below) by looking at the corresponding line graph $L(G)$ for two reasons: the spheres on $L(G)$ don't correspond to our edge spheres, and $\widetilde{L(G)}\neq L(\widetilde{G})$. (Recall that the vertices of $L(G)$ correspond to the edges of $G$, and two vertices in $L(G)$ are connected by an edge if the corresponding edges in $G$ have a vertex in common).

Finally, we look at \emph{semiregular} graphs, which are connected bipartite graphs where every edge connects a vertex of degree $p+1$ to one of degree $q+1$, so that the edge degree is constant at $p+q$. We shall only consider functions on the edges of this type of graph, as the vertex case will be treated in Section \ref{other} when we deal with bipartite graphs.
\begin{thm}\label{semihom}
Let $G$ be a semiregular simple graph with edge degree $p+q$, where $p,q\geq2$, and let $f:E\rightarrow\mathbb{R}$ be a function on its edges. Then we have for any basepoint $v_0\in \widetilde{V}$
$$\Big|M_{r,v_0}(f) - \frac{1}{|E|}\sum_{e\in E} f(e)\Big| \leq  C_G ||f||_2 \beta^r $$
Here $C_G$ is a constant depending on $G$ but independent of $v_0$, and $\beta\in[(pq)^{-1/4},1)$.
\end{thm}

Note that Theorem \ref{semihom} only deals with bipartite graphs, whereas Theorem \ref{reg e} holds for any graph, so it is not a special case of Theorem \ref{semihom}. We need $p,q\geq2$ in this theorem, as $p=1$ can give a non-converging function on the graph - see the example in Figure 2 on $K_{2,3}$. Taking the top left vertex as $v_0$, we find the spherical average takes the values $1,-1,-1,1,1,-1,\ldots$ for increasing $r$, and doesn't converge.

\begin{figure}[ht]
\centering
\includegraphics[height=3cm]{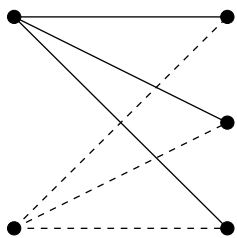}
\caption{Function on the edges of $K_{2,3}$: dashed edges have function value $-1$, others have function value $1$.}
\end{figure}

Before we turn to the proofs of the above theorems and discuss further results, let us briefly explain how spherical means are related to non-backtracking random walks (NBRW). This is a topic which has been investigated thoroughly in recent years and is the subject of active current research. 

NBRW have been studied in the context of cogrowth, although there is a subtle difference between the way we calculate our average and the method used for NBRW. In the latter case, one calculates the probability with which a non-backtracking path in $G$ with starting point $v\in V$ and length $n$ will terminate at a vertex $w$; this is calculated for all $w$, and the resulting probabilities can be used as weights in calculating the average of a function $f$ on the vertices. However in our case all non-backtracking paths are counted with the same weight, which in the case of non-regular graphs gives a different average.

The concept of cogrowth on graphs was first introduced in the context of groups and their Cayley graphs. This was studied in the early 1980s by Grigorchuk \cite{grig}, Cohen \cite{cohen} and Woess \cite{woess}. In the 1990s the application of cogrowth was extended to arbitrary graphs, see for example Northshield \cite{north3} or Bartholdi \cite{barth}. We define the growth of the tree $\widetilde{G}$ by $$\mathrm{gr}(\widetilde{G})= \limsup_{r\rightarrow\infty} \big| S_r(\tilde v) \big|^{1/r},$$ and the cogrowth of the graph $G$ by $$\mathrm{cogr}(G)= \limsup_{r\rightarrow\infty} \big| S_r(\tilde v) \cap \pi^{-1}(v) \big|^{1/r},$$ both of which are independent of $\tilde v \in \widetilde{V}$, where $\pi(\tilde v)=v\in V$. Then the cogrowth constant is $\eta = \frac{\ln \mathrm{cogr}(G)}{\ln \mathrm{gr}(\widetilde{G})}$.

More recently, Ortner and Woess \cite{ow} generalised the definition of cogrowth and used it to study NBRW. They set $$\mathrm{cog}_r^\nu(v,w)=\nu_{\tilde{v},r} (\pi^{-1}(v))$$ where $\nu=(\nu_{\tilde{v},r})_{\tilde{v}\in\widetilde{V},r\geq0}$ is a sequence of probability measures concentrated on the sphere $S_r(\tilde{v})$, subject to some regularity conditions. Choosing particular measures we obtain cogrowth or NBRW probabilities, which coincide in the case of a regular graph. In this case the results in their paper prove convergence of our spherical mean, however without giving information on the convergence rate, in contrast to the methods used in this article. In fact although various authors have investigated the spectrum of the non-backtracking adjacency matrix, which is related to the poles of the Ihara zeta function (see e.g.~\cite{st}, \cite{afh}), they have not related this to the convergence rate of the spherical mean as defined here.

The convergence of a \emph{simple} random walk (which allows backtracking) on a finite $(q+1)$-regular graph $G$ is discussed in \cite{hlw}, where the authors find that the probability distribution of the random walk converges to the uniform distribution by a factor $\alpha<1$ at each step, where $(q+1)\alpha$ is the modulus of the largest nontrivial eigenvalue of the vertex adjacency matrix of $G$. They also find the lowest value of $\alpha$ appears iff $G$ is Ramanujan. However this random walk allows backtracking, so we cannot use this result in the case of spherical means, as the paths that define the vertices of $\widetilde{G}$ do not allow backtracking.

\section{Proof of Theorem \ref{reg v}}\label{pf reg v}

We start by discussing Theorem \ref{reg v}, namely the case of a function on the vertices of a connected non-bipartite regular graph $G$ of degree $d(v)=q+1\geq3$. Let $f:V \rightarrow \mathbb{R}$ be a function on the vertices of a graph $G$. Then we define the Laplacian of $f$ at $v\in V$ as $$\mathcal{L}_G f(v)=\frac{1}{d(v)}\sum_{w\sim v} f(w)$$ where $w\sim v$ means that the vertex $w$ is adjacent to $v$. It is related to the adjacency matrix $A_G$ of the $(q+1)$-regular graph via $\mathcal{L}_G=\frac{1}{q+1} A_G$, which means it is a real symmetric operator with eigenvalues $\mu$ satisfying $-1\leq\mu\leq1$. The eigenvalue $-1$ occurs iff $G$ is bipartite, and we have excluded this case from the theorem as the result will not hold for such a graph (precisely due to this eigenvalue). The simple eigenvalue $1$ is associated to the constant eigenfunction, so for all nonconstant eigenfunctions we now have $|\mu|<1$. We will show that the spherical mean converges to the graph average, and then use the proof to calculate the convergence rate.

There is an orthonormal basis of $|V|$ eigenfunctions $\varphi_i$ of the Laplacian, that is, $\mathcal{L}_G(\varphi_i)=\mu_i\varphi_i$ and $(\varphi_i,\varphi_j)= \sum_{v\in V}\varphi_i(v)\varphi_j(v)=0$ for $i\neq j$. We start by proving the convergence of spherical means for these eigenfunctions. Let $\varphi_0$ be the constant eigenfunction, and note that here the spherical mean $M_r$ clearly equals the graph average for all $r$. Since $\sum_{v\in V}\varphi_i(v)=\sqrt{|V|}(\varphi_0,\varphi_i)=0$ for $i\neq0$, we are left to prove that $M_r(\varphi_i)\rightarrow0$ for $i=1,2,\ldots,|V|-1$ with the required speed of convergence. 

For each eigenfunction $\varphi_i\neq\varphi_0$ on $G$ let $\widetilde\varphi_i$ be its lift onto the universal covering tree $\widetilde{G}$, where it is also an eigenfunction of the Laplacian, and define the \emph{radialisation} of $\widetilde\varphi_i$ with respect to $v_0\in V$ as
$$F_i(v)= \frac{1}{|S_n(v_0)|} \sum_{w\in S_n(v_0)} \widetilde\varphi_i(w)$$
where $n=\delta(v_0,v)$. $F_i$ is also an eigenfunction on $\widetilde{G}$ with the same eigenvalue $\mu_i$ as $\varphi_i$, and as it only depends on $n$ we shall denote it $F_i(n)$ for all $v$ such that $\delta(v_0,v)=n$. Observe that $F_i(n) = M_n(\varphi_i)$. Using $\mathcal{L}_{\widetilde{G}}(F_i)=\mu_i F_i$ we obtain a recursion relation
$$F_i(n+2) - \frac{q+1}{q} \mu_i F_i(n+1) + \frac{1}{q} F_i(n) =0$$
the solution to which is given by $F_i(n) = u_+ \alpha_+^n + u_- \alpha_-^n$, where
$$\alpha_\pm = \frac{q+1}{2q}\mu_i \pm \frac{1}{2q}\sqrt{D} $$
for $D=(q+1)^2\mu_i^2-4q \neq0$, and $u_\pm$ are some constants derived from the initial values $F_i(0)$, $F_i(1)$, which in turn are determined by the eigenfunction $\varphi_i$ and the chosen basepoint $v_0$. For $D=0$, $\alpha_\pm = \alpha$ and $F_i(n)= u_1 \alpha^n + u_2 n \alpha^n $ for some constants $u_1,u_2\in\mathbb{C}$. It then just remains to check that $| \alpha_\pm | <1$ and $|n\alpha^n|\rightarrow0$ to show $\lim_{n\rightarrow\infty} F_i(n)=0$ for $i\neq0$.

For the calculation of the convergence rate we distinguish three cases:

\textbf{Case $D<0$} ($|\mu_i|<\frac{2\sqrt{q}}{q+1}$): We find $|\alpha_\pm| =\frac{1}{\sqrt{q}}$ and
$$ |F_i(n)| \leq (|u_+|+|u_-|) \big( \frac{1}{\sqrt{q}} \big)^n \leq C_i q^{-n/2} $$
for some constant $C_i>0$ which originally depends on $u_+$ and $u_-$, that is $\varphi_i$ and $v_0$. Now there are only finitely many values of $u_\pm$, because $F_i$ is the same for any lift $v_0$ of one of the \emph{finitely many} $v\in V$. Therefore we can choose $C_i$ large enough so that it is independent of $v_0$. 

\textbf{Case $D=0$} ($|\mu_i|=\frac{2\sqrt{q}}{q+1}$): Here we have
$$|F_i(n)| \leq \big( |u_1|+|u_2|n \big) q^{-n/2} \leq C'_i \cdot (n+1) \cdot q^{-n/2}$$ 
for some $C'_i>0$. Choosing $\beta_i=q^{-1/2+\varepsilon}$ for some $\varepsilon>0$ and adjusting the constant $C_i(\varepsilon)$ appropriately, we obtain
$$|F_i(n)| \leq C_i(\varepsilon) \beta_i^n$$
for $C_i(\varepsilon)>0$ independent of $v_0$. 

\textbf{Case $D>0$} ($\frac{2\sqrt{q}}{q+1}<|\mu|<1$): We find $\alpha_\pm$ are both real, and at least one of them satisfies $\frac{1}{\sqrt{q}} < |\alpha| <1$. Let $\beta_i = \text{max } \{|\alpha_+(\mu_i)|, |\alpha_-(\mu_i)|\}$, then $\frac{1}{\sqrt{q}}<\beta_i<1$ and we have
$$|F_i(n)| \leq C_i \beta_i^n$$
for some $C_i>0$ independent of $v_0$. 

A general function $f:V\rightarrow\mathbb{R}$ can be written as $f=\sum_{i=0}^{|V|-1} a_i\varphi_i$ and we obtain 
$$ \big| M_r(f) - \frac{1}{|V|} \sum_{v\in V} f(v) \big| \leq \Big| \sum_{i=1}^{|V|-1} a_i F_i(r) \Big| \leq \big( \sum_{i=1}^{|V|-1} |a_i| C_i \big) \beta_{\text{max}}^r $$
Here $\beta_{\text{max}}$ is the convergence rate obtained from the eigenvalue $\mu_{\text{max}}\neq1$ of largest modulus as then $\beta_1 \geq \beta_i \ \forall \ i\neq0$ (so the larger the spectral gap of $G$, the smaller $\beta_{\text{max}}$). Applying Cauchy-Schwarz, we obtain
$$\big|M_r(f) - \frac{1}{|V|} \sum_{v\in V} f(v) \big| \leq C_G \sqrt{\sum_{i=1}^{|V|-1} |a_i|^2 } \ \beta_{\text{max}}^r = C_G \ ||f||_2 \ \beta_{\text{max}}^r $$
where $C_G=\sqrt{|V|-2} \cdot \max_i C_i$. Note that this convergence is independent of the basepoint $v_0$, and that for Ramanujan graphs we obtain $\beta_{\text{max}}=q^{-1/2}$ (or $q^{-1/2+\varepsilon}$) as all their eigenvalues give $D\leq 0$.

\section{Proof of Theorem \ref{reg e}}\label{pf reg e}

Theorem \ref{reg e} concerns functions on the edges of a regular graph $G$, and the method of proof follows that of the vertex case apart from a small deviation towards the end. We no longer allow the graph to have loops or multiple edges, and require $d'(e)=2q \ \forall \ e \in E$. Let $g:E\rightarrow\mathbb{R}$ be a function on the edges of a graph $G$. Then the (edge) Laplacian of $g$ at $e\in E$ is defined as $$\mathcal{L}_G'g(e)=\frac{1}{d'(e)}\sum_{a\sim e} g(a)$$ 
where $a\sim e$ means that the edge $a$ is adjacent to $e$ in the sense that it has a vertex in common with $e$ (note that this is equivalent to taking the vertex Laplacian on the line graph $L(G)$ of $G$). We find that here the range of eigenvalues of the edge Laplacian is smaller, namely $-1/q \leq \mu_i \leq 1$, as for any line graph the eigenvalues of the adjacency matrix are no less than $-2$, see \cite{doob}. The recurrence relation now looks as follows
$$F_i(n+1) + \frac{q-1-2\mu_i q}{q} F_i(n) + \frac{1}{q} F_i(n-1) =0$$
Once again we want to show for $-\frac{1}{q} \leq \mu_i <1$ that $\lim_{n\rightarrow\infty} F_i(n)=0$. For $D=(q-1-2\mu_i q)^2-4q \neq 0$ we find again that $F_i(n) = u_+ \alpha_+^n + u_- \alpha_-^n$, where this time
$$\alpha_\pm = \alpha_\pm(\mu_i) = \mu_i-\frac{q-1}{2q} \pm \frac{1}{2q}\sqrt{D} $$
and for $D=0$ we have $F_i(n)= u_1 \alpha^n + u_2 n \alpha^n $. 

For $D\leq0$ the proof now follows that of the vertex case. Note that $D>0$ for $\mu\in [-1/q, \mu_1) \cup (\mu_2,1]=I$, where $\mu_1= \frac{q-1}{2q}-q^{-1/2}$ and $\mu_2= \frac{q-1}{2q}+q^{-1/2}$. As functions of $\mu$, $\alpha_\pm$ are both monotone on $[-1/q,\mu_1)$ and $(\mu_2,1]$, because $\frac{\partial}{\partial\mu} \alpha_\pm$ doesn't change sign in either interval. Calculating $|\alpha_\pm|$ for boundary values of $I$ gives $|\alpha_\pm|<1 \ \forall \ \mu \in I$, except $\alpha_+(1)=1$ (which corresponds to the constant funtion) and $|\alpha_-(-1/q)|=1$. 

However Theorem 3 in \cite{doob} states that any eigenfunction on the edges of a graph with eigenvalue $-2$ (which corresponds to $\mu=\frac{-1}{q}$) must have $\sum_{e\ni v_0} f(e)=0 \ \forall \ v_0\in V$, where $e\ni v_0$ means that $v_0$ is either the initial or terminal vertex of the edge $e$. But this is equivalent to saying that $F_i(0)=0 \ \forall \ v_0 \in V$. Observing that $2q\mu F_i(0) = qF_i(1)+qF_i(0)$ we find $F_i(n)\equiv 0 \ \forall \ n\in\mathbb{N}$, which means the spherical average of the corresponding $\varphi_i$ converges to zero as required. This completes the proof of the fact that the spherical mean of functions on the edges of $G$ converge to the graph average. To find the convergence rates we work completely analogously to the vertex case in Theorem \ref{reg v}.

\section{Proof of Theorem \ref{semihom}}\label{pf semihom}

Finally we discuss our third theorem, the case of functions on the edges of a simple semiregular graph with edge degree $p+q$, where we require that $p,q\geq2$. As for Theorem \ref{reg e}, we reduce the problem to looking at the radialisation of nonconstant eigenfunctions of the edge Laplacian, which now has eigenvalues $\frac{-2}{p+q} \leq \mu \leq 1$. Recall that edge circles are centred on a vertex, and we shall assume this vertex has degree $p+1$.

Because $G$ is semiregular, there is a more complicated recursion formula for the radialised eigenfunction $F$ with eigenvalue $\mu$ on the edges of the universal covering tree $\widetilde{G}$, given by
\begin{displaymath}
\binom{F(2k+1)}{F(2k)} = A \cdot \binom{F(2k-1)}{F(2k-2)} \qquad \text{where}
\end{displaymath}
\begin{displaymath}
A= \left( \begin{array} {ccc}
\frac{\big(p-1-\mu(p+q)\big)\big(q-1-\mu(p+q)\big)-p}{pq} & & \frac{p-1-\mu(p+q)}{pq} \\
 & & \\
-\frac{q-1-\mu(p+q)}{p} & & -\frac{1}{p} \end{array} \right)
\end{displaymath}
Hence $\binom{F(2k+1)}{F(2k)} = A^k \cdot \binom{F(1)}{F(0)}$. The convergence properties of the spherical average are now determined by the eigenvalues of the matrix $A$, which are
$$t_\pm= t_\pm(\mu)=\frac{\big(p-1-\mu(p+q)\big) \big(q-1-\mu(p+q)\big) -p-q \pm \sqrt{D(\mu)}}{2pq}$$
$$\text{where} \quad D = D(\mu)= \Big(\big(p-1-\mu(p+q)\big) \big(q-1-\mu(p+q)\big) -p-q\Big)^2 -4pq$$

It turns out below that the convergence of the spherical mean can only fail if we have $\mu_i$ such that $|t_+(\mu_i)|\geq1$ or $|t_-(\mu_i)|\geq1$ (see formulas (\ref{1}) and (\ref{2}) below). Therefore we investigate $|t_\pm(\mu_i)|$ for all possible $\mu_i$, and we distinguish the cases $D(\mu)\leq0$ and $D(\mu)>0$.

For $D\leq 0$ we have $|t_\pm|=\frac{1}{\sqrt{pq}}$ which means the spherical average converges as required. When $D>0$ we look at the various regions of $\mu$ for which $D>0$ separately. Note $D=0$ for 
$$\mu_{\pm\pm}=\frac{p+q-2\pm \sqrt{(p-q)^2 +4(\sqrt{p} \pm \sqrt{q})^2}}{2(p+q)}$$
and $D>0$ in the intervals $I_1=[\frac{-2}{p+q},\mu_{-+})$, $I_2=(\mu_{--},\mu_{+-})$ and $I_3(\mu_{++},1]$, where the subscripts $+$ and $-$ refer to the choices of $\pm$ in $\mu_{\pm\pm}$ in order of appearence. Solving $\frac{\partial}{\partial \mu} t_\pm=0$ gives $\mu'=\frac{p+q-2}{2(p+q)}$ as only solution for $\mu \in I_1 \cup I_2 \cup I_3$. Calculating $|t_+|$ at all boundary values of $I_i$ and at $\mu'$, we find $|t_+|\leq1$ with equality iff $\mu=1$ (constant eigenfunction) or $\mu=\frac{-2}{p+q}$, in which case we use Theorem 3 in \cite{doob} again to find $F(n)=0 \ \forall \ n\in \mathbb{N}$ as in Section \ref{pf reg e}.

The second eigenvalue $t_-$ is more problematic. If we have $p,q$ such that the root $\sqrt{(p-q)^2 -4(p-1)(q-1)}$ is imaginary, we obtain $|t_-(\mu)|<1$ for all $\mu$; but when $p$ and $q$ are such that it is real, we choose
$$\mu_\pm= \frac{p+q-2 \pm \sqrt{(p-q)^2 -4(p-1)(q-1)} }{2(p+q)}$$
and derive that $t_-(1)=t_-(\frac{-2}{p+q})=1$ (solved as before), $|t_-(\mu)|<1$ for $\mu\in (\frac{-2}{p+q},\mu_-) \cup (\mu_+,1)$, $t_-(\mu_\pm)=-1$, and $t_-(\mu)<-1$ for $\mu\in(\mu_-,\mu_+)$. However, as we will show in Lemma \ref{gap} below, there can be no eigenvalues $\mu$ in the interval $(\frac{p-1}{p+q},\frac{q-1}{p+q})$, where we assume that $p\leq q$ (if $p>q$ use the interval $(\frac{q-1}{p+q},\frac{p-1}{p+q})$ instead). As $[\mu_-,\mu_+] \subset \big( \frac{p-1}{p+q}, \frac{q-1}{p+q} \big)$, we find that $|t_-|<1$ for all eigenvalues $\mu$ that occur, which proves convergence of the spherical mean.

As for the convergence rate, we first assume that $D\neq0$, and let $\underline u_1$, $\underline u_2$ be a basis of eigenvectors of $A$ corresponding to the eigenvalues $t_+$, $t_-$ respectively. Writing the initial vector $\binom{F(1)}{F(0)} = a_1 \underline u_1 + a_2 \underline u_2$ we find
\begin{equation} \label{1}
\binom{F(2k+1)}{F(2k)} = A^k \binom{F(1)}{F(0)} = a_1 t_+^k \underline u_1 + a_2 t_-^k \underline u_2
\end{equation}
For $D<0$ we now use the fact that $|t_+| = |t_-| = \frac{1}{\sqrt{pq}}$ to find
$$|F(2k+j)| \leq C_j \big( \frac{1}{\sqrt{pq}} \big)^k  \qquad \text{for } j=0,1$$
with suitable constants $C_0$, $C_1$ both depending only on $F(0)$ and $F(1)$, hence
$$|F(n)| \leq C (pq)^{-\frac{n}{4}}$$
for some $C>0$ depending on $F(0)$ and $F(1)$. 

When $D>0$ the convergence will depend on the eigenvalue of $A$ with largest absolute value. Letting $\beta=\text{max}\{|t_+|,|t_-|\}$ and using the same methods as before we find
$$|F(n)| \leq C \beta^{\frac{n}{2}}$$
for some $C>0$ depending on $F(0)$ and $F(1)$, and $\frac{1}{\sqrt{pq}} < \beta<1$. 

Now in the case that $D=0$, $A$ has an eigenvalue $t=\frac{1}{\sqrt{pq}}$ or $\frac{-1}{\sqrt{pq}}$ of algebraic multiplicity two; choosing a Jordan base $\underline u_1$, $\underline u_2$ and constants $a_1$, $a_2$ appropriately such that $\binom{F(1)}{F(0)} = a_1 \underline u_1 + a_2 \underline u_2$ we derive
\begin{equation} \label{2}
\binom{F(2k+1)}{F(2k)} = A^k \binom{F(1)}{F(0)} = (a_1 t^k + a_2 k t^{k-1}) \underline u_1 + a_2 t^k \underline u_2
\end{equation}
This implies that
$$|F(n)| \leq C' \cdot (1+n) \cdot (pq)^{-\frac{n}{4}} \leq C \beta^{\frac{n}{2}}$$
for $C'>0$ depending on $F(0)$ and $F(1)$, and $\beta=(pq)^{-\frac{1}{2}+\varepsilon}$ for arbitrarily chosen $\varepsilon>0$, and appropriately adjusted $C$.

As with the previous two theorems, we write $f=\sum_{i=0}^{|V|-1}b_i\varphi_i$ and use the largest value of $\beta$ to find $$\Big|M_{r,v_0}(f)-\frac {1}{|E|}\sum_{e\in E} f(e)\Big| \leq C_G ||f||_2 \beta^r$$
where as before  $C_G>0$ large enough to provide independence of $v_0\in V$. 

To complete the proof of Theorem \ref{semihom}, it remains to prove the following lemma.
\begin{lem}\label{gap}
Let $G$ be a semiregular graph as in Theorem \ref{semihom}. Then the edge Laplacian has no eigenvalues $\mu$ such that $$\frac{p-1}{p+q} <\mu< \frac{q-1}{p+q}$$ 
\end{lem}

\begin{pf}
Let $G$ be a semiregular graph with $n_1$ vertices of degree $p+1$ and $n_2$ vertices of degree $q+1$, where $n_1\geq n_2$ and all vertices of the same degree are mutually non-adjacent. Then a theorem by Cvetkovi\'c (Theorem 1.3.18 in \cite{crs}) gives the following relation between the characteristic polynomials $P_G(x)$ and $P_{L(G)}(x)$ of $G$ and its line graph $L(G)$ respectively:
$$ P_{L(G)}(x) = 
(x+2)^{m}  \sqrt{\Bigg(\frac{-\alpha_1(x)}{\alpha_2(x)}\Bigg)^{n_1-n_2} P_G\Big(\sqrt{\alpha_1(x)\alpha_2(x)}\Big) P_G\Big(-\sqrt{\alpha_1(x)\alpha_2(x)}\Big)}$$
where $m=|E|-|V|$, $\alpha_1=x-p+1$ and $\alpha_2=x-q+1$. Recall $P_{L(G)}(\lambda)=0$ for eigenvalues $\lambda$ of the edge adjacency matrix $A_{L(G)}$, and as $\mathcal{L}'_{L(G)}=\frac{1}{p+q}A_{L(G)}$ we have
$$\mu=\frac{\lambda}{p+q}$$
so $\lambda\in[-2,p+q]$. Using the above formula for $P_{L(G)}$, we find its roots can only be $\lambda=-2$, $\lambda=p-1$, or values of $\pm\sqrt{\alpha_1(x)\alpha_2(x)}$ that are eigenvalues of the original graph. Note that $G$ has only real eigenvalues. However since $\sqrt{\alpha_1(x)\alpha_2(x)}$ is \emph{purely imaginary} for $p-1<x<q-1$, $L(G)$ cannot have eigenvalues in this region. \hfill $\square$
\end{pf}

\section{Further Results}\label{other}

In this section we discuss some related results.

Let us first revisit the case of bipartite graphs, which were excluded when we investigated functions on the vertices of a regular graph. We shall also consider semiregular graphs, with functions defined on the vertices, and we shall allow both types of graph to have multiple edges. The problem in Theorem \ref{reg v} was the fact that there exists an eigenfunction of the Laplacian on a bipartite graph with eigenvalue $-1$, namely the function with value $A$, $-A$ respectively on the respective parts of the graph. The spherical average then alternates between the two values as $n\rightarrow\infty$ and never converges. We remedy this by modifying the graph and looking at the two parts seperately in the following way. Let $V_1$ and $V_2$ be the two independent sets of vertices (i.e.~$x\nsim y\ \forall \ x,y\in V_i$). Without loss of generality, suppose $v_0 \in V_1$. Define a new graph $G'$ with vertex set $V_1$, and an edge connecting $x$ and $y$ in $G'$ for each non-backtracking path of length two between $x$ and $y$ in $G$. This gives us a regular non-bipartite graph $G'$ of degree $q(q+1)$, to which we can apply Theorem \ref{reg v}. Applying the same method to a semiregular graph gives a regular graph of degree $(p+1)q$ or $p(q+1)$, depending on our choice of $V_1$ (vertices of degree $p+1$) or $V_2$ (vertices of degree $q+1$), and again we can apply Theorem \ref{reg v}. Therefore, for a bipartite regular or semiregular graph, we have two limits of spherical averages for increasing \emph{even} radii, depending on which part of the graph contains the base vertex.

\begin{figure}[ht]
\centering
\includegraphics[height=3.5cm]{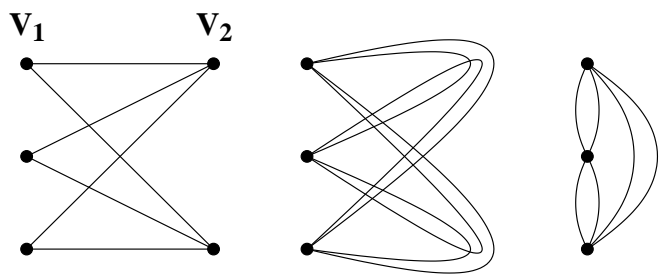}
\caption{From left to right: the transition from $G$ to $G'$}
\end{figure}

Secondly we find, as mentioned in Section \ref{intro}, that we can in fact extend the result in Theorem \ref{reg v} to increasing subsets of $\widetilde{G}$ other than spheres. 

We define a circular \emph{arc} $A_r$ on the tree as follows. Let $\vec{e}$ be a directed edge from vertex $w'$ to $w$, then $A_{r}(\vec{e})=S_{r+1}(w') \cap S_r(w)$ is based at $w$ in the direction of $\vec{e}$ (see Figure 4). Now let $X$ be a connected subgraph of $T$. We define the \emph{tube} $\mathcal{T}_r$ in $T$ of radius $r$ around $X$ as 
$$\mathcal{T}_r(X)= \{ v \in V(T) : \min_{x \in V(X)} \delta (v,x) =r \} $$
where $\delta$ is the combinatorial distance (see also Figure 4, where $X$ consists of the thick vertices and edges).

\begin{figure}[ht]
\centering
\includegraphics[height=6cm]{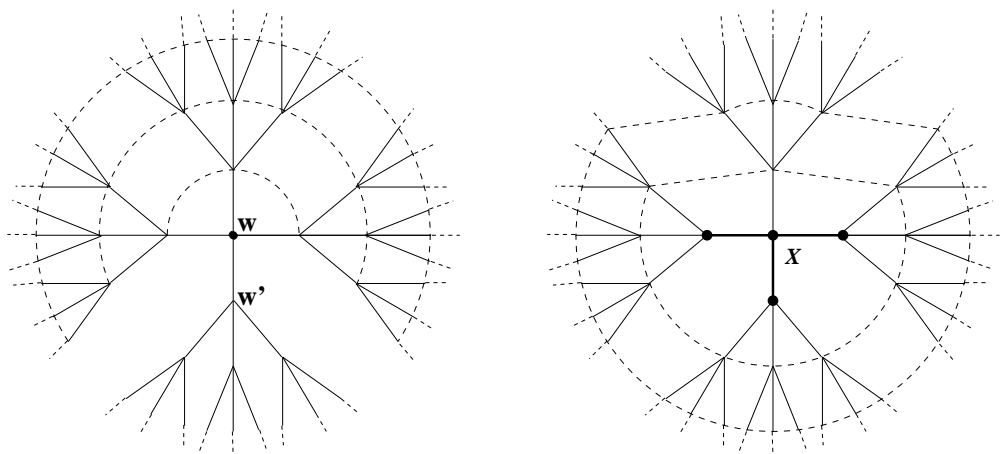}
\caption{Vertex arcs and tubes on a regular graph of degree 4}
\end{figure}

\begin{lem}\label{arcl}
 Let $G$ be a finite regular non-bipartite graph as defined in Section \ref{intro} and $f:V\rightarrow\mathbb{R}$ a function on its vertices. Then for any directed edge $\vec{e}$
$$\Big| \frac{1}{|A_{r}(\vec{e})|} \sum_{v\in A_{r}(\vec{e})} \tilde{f}(v) - \frac{1}{|V|} \sum_{v\in V} f(v) \Big| \leq C_G ||f||_2 \beta^r $$
where $C_G>0$ depends only on $G$, and $\beta \in [q^{-1/2},1)$ with lowest value $\beta=q^{-1/2}$ for Ramanujan graphs.
\end{lem}
\begin{pf}
 We use the recursion formula from Theorem \ref{reg v}, and as its convergence does not depend on the values of $F(n)$ for $n=0,1$, the same calculations give us convergence of the arc average with the same convergence rates. \hfill $\square$
\end{pf}

\begin{rmk}
 We can do the same for edge arcs, although we have to check what happens for $\mu=-1/q$: use $\sum_{e \ni v_0} f(e) =0 \ \forall \ v_0\in V$ from \cite{doob} as before to find $\tilde{f}(\vec{e}) +qF(1) =0$. Use $\mathcal{L}' F(1) = -\frac{1}{q} F(1)$ to obtain $F(n)=(-1/q)^n \tilde{f} (\vec{e})$ which clearly converges to zero as $n \rightarrow \infty$ with $\beta = 1/q$.
\end{rmk}

\begin{cor}\label{tubes}
 Let $G$ and $f$ as above. Then 
$$\Big| \frac{1}{|\mathcal{T}_r(X)|} \sum_{v\in \mathcal{T}_r(X)} \tilde{f}(v) - \frac{1}{|V|} \sum_{v\in V} f(v) \Big| \leq C ||f||_2 \beta^r $$
where $C>0$ depends on $G$ and $X$, and $\beta \in [q^{-1/2},1)$.
\end{cor}
\begin{pf}
 Tubes can be viewed as a disjoint union of several arcs (each with different $\vec{e}$), so the result follows from Lemma \ref{arcl}. \hfill $\square$
\end{pf}

Finally, we consider increasing subsets of horocycles on $\widetilde{G}$ to find a discrete analogue of a result by Furstenberg \cite{furst} on the unique ergodicity of the horocycle flow (see also chapter IV in \cite{bekka}). Recall a geodesic $\gamma$ on the tree $\widetilde{G}$ is a bi-infinite non-backtracking path, which we shall denote by its vertices $\ldots, v_{-2},v_{-1},v_0,v_1,v_2,\ldots \in \widetilde{V}$, where $v_i \sim v_{i+1}$ and $v_i\neq v_{i+2} \ \forall \ i\in\mathbb{Z}$. Recall $\delta(v,w)$ is the combinatorial distance between vertices $v$ and $w$, and define the Busemann function
$$b_{\gamma,v_k} (w) = \lim_{n\rightarrow\infty} \delta(w,v_{k+n}) -n$$
For $k\in\mathbb{Z}$ we then define the horocycle $H_k = b_{\gamma,v_0}^{-1}(k)$ (for illustration of horocycles, see also \cite{figa} Chapter I Section 9). We shall consider subsets of the horocycle $H_0$ defined by
$$\mathcal{H}_{\gamma,r}(v_0) = H_0 \cap S_r(v_r)$$
as seen in Figure 5.

\begin{figure}[ht]
\centering
\includegraphics[height=3.5cm]{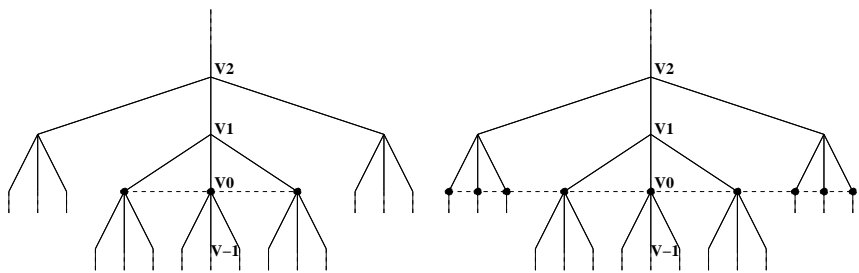}
\caption{Subsets $\mathcal{H}_{\gamma,1}(v_0)$ and $\mathcal{H}_{\gamma,2}(v_0)$ of the horocycle $H_0$}
\end{figure}

\begin{thm}
 Let $G$ and $f$ as above. Then
$$\Big| \frac{1}{|\mathcal{H}_{\gamma,r}(v_0)|} \sum_{v\in \mathcal{H}_{\gamma,r}(v_0)} \tilde{f}(v) - \frac{1}{|V|} \sum_{v\in V} f(v) \Big| \leq C_G ||f||_2 \beta^r $$
where $C_G$ depends only on $G$, and $\beta \in [q^{-1/2},1)$.
\end{thm}
\begin{pf}
 Note that we can view the subset of the horocycle as an arc
$$\mathcal{H}_{\gamma,r}(v_0) = A_r \Big(\overrightarrow{ \{v_{r+1},v_r\} }\Big) $$
so as $r\rightarrow\infty$ we have a set of increasing circular arcs, where the origin of the arc changes at each step. But the convergence for arcs in Lemma \ref{arcl} is independent of the origin of the arc, so the subsets can be viewed just as increasing circular arcs, and the theorem follows. \hfill $\square$
\end{pf}

\setlength{\parskip}{0mm}

\bibliographystyle{99}

\end{document}